\input amstex
\documentstyle{amsppt}

\magnification=1200
\loadmsbm
\loadmsam
\loadeufm
\input amssym
\UseAMSsymbols

\TagsOnRight

\hsize168 true mm
\vsize220 true mm
\voffset=20 true mm
\hoffset=0 true mm
\baselineskip 9 true mm plus0.4 true mm minus0.4 true mm

\def\l {\ell}

\topmatter
\title
From $2D $ Toda hierarchy to conformal map for domains
of Riemann sphere
\endtitle

\author 
Yu.Klimov, A.Korzh, S.Natanzon
\endauthor

\affil
Moscow State University, Independent University of Moscow\\
Moscow State University, Independent University of Moscow\\
Moscow State University, Independent University of Moscow,
ITEP
\endaffil

\endtopmatter
\document

\subhead
1. Introduction
\endsubhead

In recent works [1,2,3] was found a wonderful correlation between integrable
systems and meromorphic functions. They reduce a problem of effictivisation
of Riemann theorem about conformal maps to calculation of a string solution of
dispersionless limit of the $2D$ Toda hierarchy. In [4] was found
a recurrent formulas for coeffciens of Taylor series of the string solution.
This gives, in particular, a method for calculation of the univalent
conformal map from the until disk to an arbitrary domain, described by
its harmonic moments.

In the present paper we investigate some properties of these
formulas. In particular, we find a sufficient condition for convergence of
the Taylor series for the string solution of dispersionless limit of $2D$
Toda hierarchy.

We thank A.Marshakov and A.Zabrodin for useful discussions.

\vskip 0.4cm
\subhead
2. Complex domains and $2D$ Toda hierarchy
\endsubhead

In this section we remind some results of works [1, 2, 3].

We shall consider only domains generated by closed
analytical contours $ \gamma $ without self--intersections
on Riemann sphere $ \bar\Bbb C =\Bbb C\cup\infty $.
The analyticity of $\gamma$ is means, that $\gamma$ is
an image of $ \gamma_0 = \{w\in\Bbb C || w | = 1 \} $
by a function and this function is analytical in a
nieghborhood of $ \gamma_0 $. A curve
$ \gamma $ divides $ \Bbb C $ into an interior domain $D _ + $
and an external domain $D_-$. We shall consider, that $D _ +\ni 0 $.

Harmonic moments of $D_+$ are  the numbers 
$$ v _0 =\frac 2\pi \int_{D_+}\log |z|d^2 z, \ v_k
=\frac 1\pi\int_{D_+}z^k d^2 z\ . $$

Harmonic moments of $D_- $ are the numbers 
$$t _0=\frac 1\pi \int_{D_+}d^2z, \ t_k =
-\frac{1} {\pi k} \int_{D_-}z^{-k} d^2 z\ . $$

Let $T$ be the set of analytic curves on $\bar{\Bbb C}$,
parametrized by collections $\{t_i \}$ of its harmonic moments.
Consider on $T$ the function
$$F(t)=-\frac 1\pi\iint_{D_+(t)}\log\left\vert\frac 1z-\frac {1}{z'}
\right\vert d^2zd^2z'\ .$$
Later we consider that $F$ is a real analytical function from
$\tilde t = (t_0,t_1,\bar t_1,t_2,\bar t_2,...)\ .$

Consider the function
$$\varphi(z,\tilde t)=-\partial _ {t_0} \left (\frac 12\partial_{t_0} +
\sum _ {k\geqslant 1}\frac{z^{-k}}{k}\partial_{t_k}\right)F(\tilde t)\ ,$$
where $\partial_t=\frac{\partial}{\partial t}$ and
$z\in \bar{\Bbb C}$. Then the function
$w(z)=e^\varphi z$ is the one-sheeted meromorphic function
mapping from $D_-(t)$ to
$\{w\in\bar {\Bbb C} \left\vert|w |>1\right.\}.$ Therefore, if we
know a Taylor series of $F(\tilde t)$, we can find the functions
$p_j(t)$, giving a representation of the function $w(z)$
in the form
$$w(z)=\frac 1r z+\sum^\infty_{j=0} p_j(t) z^{-j}\ .$$

Thus, for an effectivisation of the Riemann theorem it is
sufficiently to find the Taylor series of $F(\tilde t)$.

In [1 - 3] is proved that $F(\tilde t)$ satisfy the differential
equations: 
$$(z-\xi)e^{D(z)D(\xi)F}=ze^{-\partial_{0}D(z)F}- 
\xi e^{-\partial_{0} D (\xi) F}, \tag1 $$
$$(\bar z-\xi) e^{\bar D (\bar z) \bar D (\bar\xi) F} = \bar z
e^{-\partial_{0} \bar D (\bar z) F} -
\bar\xi e^{-\partial_{0} \bar D (\bar\xi) F}, \tag2 $$
$$1-e^{D (z) \bar D (\bar\xi) F} = \frac {1} {z\bar\xi}
e^{-\partial_{0} (\partial_{0} + D(z)+\bar D(\bar\xi))F}, \tag3$$
where $\partial_k=\partial_{t_k}$,
$\bar\partial_k=\partial_{\bar t_k}$ and
$$D(z)= \sum_{k\geqslant 1} \frac {z^{-k}}{k}\partial_{k},   
\bar D(\bar z)=\sum_{k\geqslant 1}\frac{\bar z^{-k}}{k}\bar\partial_{k}\ .$$

This system of the nonlinear differential equations is well known in
mathematical physics and in theory of integrable systems as
dispersionless limit of the $2D$ Toda hierarchy [6]. 
The solution $F(\tilde t)$ satisfies some additional equation,
which appear in string theory and, therefore, $F(\tilde t)$ is called
"string solution" [7]. The string solution of dispersionless
limit of $2D$ Toda hierarchy is appeared also in
matrix models and in some other problems of mathematical
physics. Thus a description of it has an independent
interest.

\subhead
3. Taylor series for the string solution of dispersionless limit
of the $2D $ Toda \linebreak hierarchy
\endsubhead

In [ 4 ] was found a representation of $F$ in form of Taylor
series 
$$ F =\sum N(i|i_1..., i_k |\bar i_1..., \bar i_k)
t_0^i t _ {i_1} \dotsb t _ {i_k} \bar t _ {\bar i _ {\bar 1}}
\dotsb\bar t _ {\bar i _ {\bar k}}\ .$$

The formulas for $N$ are found by the following scheme.
At first, using some combina\-torial calculations, we transform
the equation (1) to an infinite system of equations
 $$\partial_{i_1}\partial_{i_2}\dotsb\partial_{i_k}F\ =$$ 
$$= \ \sum^\infty _ {m=1} \left (\sum\Sb s_1 +\dotsb +s_m =
i_1 +\dotsb +i_k \\ \ell_1 +\dotsb + \ell_m
=m+k-2 \\ s_j, \ell_j\geqslant 1\endSb
\frac {i_1\dotsb i_k} {s_1\dotsb s_m}
T _ {i_1\dotsb i_k} \pmatrix s_1\dotsb s_m \\
\ell_1\dotsb \ell_m\endpmatrix\partial_0^{\ell_1}\partial_{s_1} 
F\dotsb \partial_0 ^ {\ell_m} \partial _ {s_m} F\right). \tag4 $$

In passing we find some recurrent formulas for calculation of $T$.

Later, using the definition of $F$ as a function on the space
of analytical curves, we find, that

$\partial_0 F\left\vert _ {t_0} \right. =-t_0+t_0\ln t_0$ and
$\partial_k F\left\vert _ {t_0} \right. = 0, $ if $k > 0,$

\noindent
where here and later $\left |_{t_0}\right.$ means the restriction
of a function on a straight line
$t_1 =\bar t_1=t_2 =\bar t_2 =\dotsb=0 $.

For this formula and from the equation (3) follow, that
$$\partial_i\bar\partial_j F\left | _ {t_0} \right. =\cases 0, \
\text {if} \ i\ne j, \\ it^i_0,\ \text {if} \ i=j\ .\endcases $$

Later, using (4) and the symmetry of the equations (1) - (3)
we find, that
 $$\partial_i\bar\partial_{i_1}\dotsb\bar\partial_{i_k} 
 F\left|_{t_0}\right.=\bar\partial_i\partial_{i_1}\dotsb\partial_{i_k} 
F\left | _ {t=0} \right. =\cases 0, \ \text {if} \ i_1 +\dotsb+i_k\ne i, \\
i_1\dotsb i_k\frac {i!} {(i-k+1)!} t_0^{i-k+1},
\ \text {if}\  i_1 +\dotsb+i_k = i.\endcases$$

This condition and the equation (4) give some
recurrent formulas for coefficients $N$. As the final result we get

\proclaim {Theorem 1} In the domain of its convergence the formal series

$F\ = \ \frac {1} {2}\, t_0^2 \, \log {t_0}\ - \ \frac {3}{4}\, t_0^2+$

\vskip 0.4cm
\noindent\qquad
$ + \sum\limits_{\Sb k, \bar k, n_r, \bar n_r \geqslant 1 \\
0 < i_1 < \dots < i_k \\ 0 < \bar i_1 < \dots < \bar i _ {\bar k}
\\ i-(n_1+\dots+n_k+\bar n_1+\dots+\bar n_{\bar k})+2\geqslant 0
\endSb}
\frac {i_1 ^ {n_1} \dots i_k^{n_k}} {n_1! \dots n_k!} \ 
\frac {\bar i_1 ^ {\bar n_1} \dots \bar i_{\bar k}^{\bar n _ {\bar k}}}
{\bar n_1! \dots \bar n_{\bar k}!}
\ N_i^2 \left (\matrix i_1,&\dots,&i_k \\ n_1,&\dots,&n_k \endmatrix 
\left | \matrix \bar i_1,&\dots,&\bar i _ {\bar k} \\ \bar n_1,&\dots,& 
\bar n _ {\bar k} \endmatrix \right.\right) \times $

\vskip 0.4cm
\noindent\qquad
$ \times t_0^{i-(n_1+\dots+n_k+\bar n_1+\dots+\bar n_{\bar k}) 
+2} t_{i_1}^{n_1} \dots t_{i_k}^{n_k} \bar 
t_{\bar i_1}^{\bar n_1} \dots \bar t_{\bar i_{\bar k}}
^{\bar n_{\bar k}}$

\noindent
is the string solution
of dispersionless limit of $2D$ Toda hierarchy. In this formula
the coefficients $N^2$ are found by follow recurrent rules.

\noindent
$P_{i, j}(s_1, \dots, s_m) \ =
\ \# \ \{(i_1, \dots, i_m) \ | \ i = i_1 + \dots + i_m, \ 
1 \leqslant i_r \leqslant s_r-1 \},$ where $\# Q$ is the cardinality
of the set $Q$;

\noindent
$T_{i, j} ^1 (s_1, \dots, s_m) \ = $

$ \sum\limits_{\Sb k \geqslant 1 \\ n_1 + \dots + n_k=m \\ n_r \geqslant 1 \endSb}
\frac {1} {k  n_1! \dots n_k!} \
P_{i, j} \left (\underbrace {s_1+\dots +s_{n_1}}_{n_1},
\dots, \underbrace {s_{n_1+\dots+n_{k-1}+1}+\dots+s_{n_1+\dots+n_k}}_{n_k}
\right);$

\noindent
$T_{i_1, i_2}^2 \left (\matrix s_1,&\dots,&s_m \\ 1,&\dots,&1 \endmatrix
\right) \ = \ T_{i_1, i_2}^1 (s_1, \dots, s_m);$

\noindent
$T_{i_1, \dots, i_k}^2 \left (\matrix s_1,&\dots,&s_m \\ l_1,&\dots,&l_m  
\endmatrix\right) \ = $

\noindent
$ = \ \sum\limits_{\Sb 1 \leqslant i \leqslant j \leqslant m \\ s, l \geqslant 1 \endSb}
l \ T_{s, i_k}^1(s_i, \dots,  s_j)
\ T_{i_1, \dots, i_{k-1}}^2 \left (\matrix s_1,&\dots,&s_{i-1},&s, 
 &s_{j+1},&\dots,&s_m \\ l_1,&\dots,&l_{i-1},&l,&l_{j+1},&\dots,&l_m  
\endmatrix\right)$\ ,

\noindent
where
$s = s_i + \dots + s_j - i_k, \ l = (l_i-1) + \dots + (l_j-1)\ ;$

\noindent\qquad
$S_{\bar i_1, \dots, \bar i_{\bar k}} \left (\matrix s_1,&\dots,&s_m  
\\ l_1,&\dots,&l_m \endmatrix\right) \ = $

\vskip 0.4cm
\noindent
$=\ \sum\limits_{\Sb \{\bar i_1^1, \dots, \bar i_1^{n_1} \} \sqcup \dots 
\sqcup \{\bar i_m^1, \dots, \bar i_m^{n_m} \} = 
\{\bar i_1, \dots, \bar i_{\bar k} \} \\ \bar i_r^1 + \dots + 
\bar i_r ^ {n_r} = s_r \\s_r-n_r-\ell_r+1\geqslant 0 \endSb} 
\frac {(s_1-1)!} {(s_1-n_1-l_1+1)! (l_1-1)!} \times \dots \times\
\frac {(s_m-1)!} {(s_m-n_m-l_m+1)! (l_m-1)!}\ ;$

\bigskip\bigskip
\noindent\qquad
$N_i^1 (i_1, \dots, i_k | \bar i_1, \dots, \bar i_{\bar k}) \ 
= \ 0, \ \text { if}\ i \ne i_1+\dots+i_k \ \text {or}\ i \ne \bar i_1
+ \dots + \bar i _ {\bar k}\ ;$

\vskip 0.4cm
\noindent
\text {in the other cases}

\vskip 0.4cm
\noindent
$N_i^1 (i | \bar i_1, \dots,  \bar i_{\bar k}) \ = 
\ \frac {(i-1)!} {(i-\bar k+1)!}\ ;$

\vskip 0.4cm
\noindent
$N_i^1 (i_1, \dots,  i_k  |  \bar i) \ = \ \frac {(i-1)!} {(i-k+1)!}\ ;$

\vskip 0.4cm
\noindent
$N_i^1 (i_1, \dots, i_k |  \bar i_1, \dots, \bar i _ {\bar k}) \ =$

\vskip 0.4cm
\noindent
$ = \ \sum\limits_{\Sb m \geqslant 1 \\ s_1 + \dots + s_m=i_1 + \dots + i_k \\ 
l_1 + \dots + l_m = m+k-2 \\ s_r, l_r \geqslant 1 \endSb}
(-1) ^ {m+1} \ S _ {\bar i_1, \dots, \bar i _ {\bar k}} \left (\matrix 
s_1,&\dots,&s_m \\ l_1,&\dots,&l_m \endmatrix\right)\ $
$T_{i_1, \dots, i_k} ^2 \left (\matrix s_1,&\dots,&s_m \\
 l_1,&\dots,&l_m \endmatrix\right)\ ,$ if $\ k, \bar k>1\ ;$

\vskip 0.4cm
\noindent
$N_i^2 \left (\matrix i_1,&\dots,&i_k \\ n_1,&\dots,&n_k \endmatrix 
\left | \matrix \bar i_1,&\dots,&\bar i _ {\bar k} \\ 
\bar n_1,&\dots,&\bar n _ {\bar k} \endmatrix \right.\right) \ = $

\vskip 0.2cm
\noindent\qquad
$=\ N_i^1 \left (\underbrace {i_1,  \dots, i_1} _ {n_1}, \dots, 
\underbrace {i_k, \dots, i_k}_ {n_k} \left | \underbrace {\bar i_1, 
\dots, \bar i_1}_{\bar n_1}, \dots, \underbrace {\bar i_{\bar k}, 
\dots, \bar i_{\bar k}}_{\bar n{\bar k}} \right.\right)\ . $
$\square$
\endproclaim

For $t_3 =\bar t_3=t_4 =\bar t_4 =\dotsb=0 $ this theorem goes to
the formula from [2]
$$ F =-\frac 34 t^2_0+\frac 12 t^2_0\ln\left(\frac{t_0}{1-4|t_2 |^2}
\right)+\frac{t_0}{1-4|t_2|^2}(|t_1|^2+t^2_1\bar t_2+\bar t^2_1 t_2)\ .$$

First two authors construct a computer program, calculating
any coefficient $N^2_i (...) $. Cal\-culations by this
program lead to hypothesis, that all coefficients
$N^2_i (...) $ are nonnegative.

\subhead
4. Some properties of the coefficients for the series $F$
\endsubhead

The combinatorial coefficients $N^2_i (...)$ have some 
remarkable properties. For example,

\vskip 0.3cm
\proclaim {Theorem 2}
$N^2_i\left (\matrix i_1\dotsb i_k \\ n_1\dotsb n_k\endmatrix
\left |\matrix 1 \\ \bar n_1 \endmatrix \right.\right)=
\cases(i-1)!,\ \text {if} \
k=n_1=1, \ i=i_1 =\bar n_1, \\ 0\ \text {in the other cases}.
\endcases $
\endproclaim

\demo {Proof} According to our definition

$$S_{\underbrace{1,  \dots, 1}_{\bar k}} \left(\matrix s_1,&\dots,&s_m
\\ l_1,&\dots,&l_m \endmatrix\right) \ =$$
$$= \ \sum\limits_{\Sb \{\bar i_1^1, \dots, \bar i_1^{n_1}\}
\sqcup \dots \sqcup \{\bar i_m^1, \dots, \bar i_m^{n_m}\} =\\=
\{\underbrace{1, \dots, 1}\} \\ \bar i_r^1 + \dots
+ \bar i_r^{n_r} = s_r \\s_r-n_r-\ell_r+1\geqslant 0 \endSb}
\frac{(s_1-1)!}{(s_1-n_1-l_1+1)!(l_1-1)!} \times \dots \times$$
$$\times\dots\times\frac{(s_m-1)!}{(s_m-n_m-l_m+1)!(l_m-1)!} \ =
\delta_{\ell_1,1}\dots\delta_{\ell_m,1}\frac{\bar k!}{s_1
\dots s_m}.$$

Thus if $k>2$, then

$$N_i^1 (i_1, \dots, i_k  |  1, \dots, 1) \ =$$
$$= \ \sum\limits_{\Sb m \geqslant 1 \\ s_1 + \dots + s_m =
i_1 + \dots + i_k \\ l_1 + \dots + l_m = m+k-2 \\ s_r, l_r \geqslant 1
\endSb}
(-1)^{m+1} \ S_{1, \dots, 1} \left(\matrix s_1,&\dots,&s_m
\\ l_1,&\dots,&l_m \endmatrix\right)\times$$
$$\times \ T_{i_1, \dots, i_k}^2 \left(\matrix s_1,&\dots,&s_m \\
l_1,&\dots,&l_m \endmatrix\right)=0\ .$$
Let now $k=2$ ($i_1, i_2 \geqslant 1$). Then
$$N_i^1 (i_1, i_2  | \underbrace{ 1, \dots, 1}_{\bar k}) \ =$$
$$= \ \sum\limits_{\Sb m \geqslant 1 \\ s_1 + \dots + s_m = i_1 + i_2
\\ s_r \geqslant 1 \endSb}
(-1)^{m+1} \ S_{\underbrace{ 1, \dots, 1}_{\bar k}}
\left(\matrix s_1,&\dots,&s_m \\ 1,&\dots,&1 \endmatrix\right)\times$$
$$\times\ T_{i_1, i_2}^2 \left(\matrix s_1,&\dots,&s_m
\\ 1,&\dots,&1 \endmatrix\right) \ =$$
$$= \ \sum\limits_{\Sb m \geqslant 1 \\ s_1 + \dots + s_m =
i_1 + i_2 \\ s_r \geqslant 1 \endSb}
(-1)^{m+1}  \frac{\bar k!}{s_1  \dots  s_m}
\ T_{i_1, i_2}^2 \left(\matrix s_1,&\dots,&s_m \\ 1,&\dots,&1
\endmatrix\right) \ =$$
$$= \ \sum\limits_{\Sb m \geqslant 1 \\ s_1 + \dots + s_m = i_1 + i_2
\\ s_r \geqslant 1 \endSb}
(-1)^{m+1}   \frac{\bar k!}{s_1  \dots  s_m}
\ T_{i_1, i_2}^1 (s_1,  \dots,  s_m) \ =$$
$$= \ \sum\limits_{\Sb m \geqslant 1 \\ s_1 + \dots + s_m = i_1 + i_2
\\ s_r \geqslant 1 \endSb}
(-1)^{m+1}  \frac{\bar k!}{s_1  \dots  s_m}  
 \sum\limits_{\Sb k \geqslant 1 \\ n_1 + \dots + n_k=m
\\ n_r \geqslant 1 \endSb}
\frac {1}{k  n_1!  \dots  n_k!} \times$$
$$\times \ P_{i, j} \left(\underbrace{s_1 + \dots + s_{n_1}}_{n_1},
\dots, \underbrace{s_{n_1 + \dots + n_{k-1} + 1} + \dots +
s_{n_1 + \dots + n_k}}_{n_k}\right) \ =$$

$$= \ \sum\limits_{\Sb m \geqslant 1 \\ s_1 + \dots + s_m = i_1 + i_2
\\ k \geqslant 1 \\ n_1 + \dots + n_k=m \endSb}
\frac {(-1)^{m+1}}{k  n_1!  \dots  n_k!}  \frac{\bar k!}{s_1
 \dots  s_m}\times$$
$$\times\ P_{i, j} \left(\underbrace{s_1 + \dots + s_{n_1}}_{n_1},
\dots, \underbrace{s_{n_1 + \dots + n_{k-1} + 1} + \dots +
s_{n_1 + \dots + n_k}}_{n_k}\right) \ =$$
$$= \ \sum\limits_{\Sb k \geqslant 1 \\ \tilde s_1 + \dots + \tilde s_k
= i_1 + i_2 \endSb}
\ P_{i, j} (\tilde s_1, \dots,  \tilde s_k)\times$$
$$\times\ \sum\limits_{\Sb n_r \geqslant 1,\ n_1+\dots+n_k = m
\\ \underbrace{s_{n_1 + \dots + n_{r-1} + 1} + \dots +
s_{n_1 + \dots + n_r}}_{n_r} = \tilde s_r \endSb}
\frac {(-1)^{m+1}}{k  n_1!  \dots  n_k!} 
\frac{\bar k!}{s_1  \dots  s_m} \ =$$
$$= \ \sum\limits_{\Sb k \geqslant 1 \\ \tilde s_1 + \dots +
\tilde s_k = i_1 + i_2 \endSb}
-\frac{\bar k!}{k}\ P_{i, j} (\tilde s_1, \dots,  \tilde s_k)\times$$
$$\times\prod\limits_{1 \leqslant r \leqslant k} \sum\limits_{\Sb n_r \geqslant 1
\\ s_1 + \dots + s_{n_r} = \tilde s_r \endSb}
\frac {(-1)^{n_r}}{n_r!  s_1  \dots  s_{n_r}}\ .$$

In addition if $s>1$, then 

$$\sum\limits_{\Sb n \geqslant 1 \\ s_1 + \dots + s_{n} =
s \endSb} \frac {(-1)^{n}}{n!  s_1  \dots  s_{n}} \ =
\ \frac{1}{s!}\frac{\partial^s}{\partial x^s}\
\sum\limits_{n \geqslant 1} \frac {(-1)^{n}}{n!} \left.
\left(x+\frac{x^2}{2}+\frac{x^3}{3}+\dots \right)^n \right|_{x=0}\ =$$
$$=\ \frac{1}{s!}\frac{\partial^s}{\partial x^s}\
\sum\limits_{n \geqslant 1} \frac {(-1)^{n}}{n!} \left. (-\log(1-x))^n
\right|_{x=0} \ = $$
$$= \ \frac{1}{s!}\frac{\partial^s}{\partial x^s}\
\left. \sum\limits_{n \geqslant 1} \frac {(\log(1-x))^n}{n!}
\right|_{x-0} \ =
\ \frac{1}{s!}\frac{\partial^s}{\partial x^s}\
\left. ( \exp(\log(1-x))-1) \right|_{x=0} \ = \ 0.$$
Thus 
$$N_i^1 (i_1, i_2  | \underbrace{ 1, \dots, 1}_{\bar k}) \ =
\ \sum\limits_{\Sb k \geqslant 1 \\ \tilde s_1 + \dots + \tilde s_k
= i_1 + i_2 \endSb}
-\frac{\bar k!}{k}\ P_{i, j} (\tilde s_1, \dots,  \tilde s_k)\times$$
$$\times\ \prod\limits_{1 \leqslant r \leqslant k} \sum\limits_{\Sb n_r \geqslant 1
\\ s_1 + \dots + s_{n_r} = \tilde s_r \endSb}
\frac {(-1)^{n_r}}{n_r!  s_1  \dots  s_{n_r}}=0\ ,$$
because
$P_{i, j} (1, \dots,  1) \ = 0.$

If $k=1$, then from our definion it follows that

$$N^2_i\left (\matrix i_1 \\ n_1\endmatrix
\left |\matrix 1 \\ \bar n_1 \endmatrix \right.\right)
=\cases(i-1)!,\ \text {if} \
n_1=1, \ i=i_1 =\bar n_1, \\ 0\ \text {in the other cases}.
\endcases\ .\ \qed$$
\enddemo

\subhead
5. Convergence conditions for the Taylor series
\endsubhead

The recurrent formulas for coefficients of the Taylor series $F$ give
possibility to estimate the coefficients and to find sufficient convergence
conditions for $F$.

\proclaim
{Theorem 3} Let $\tilde t = (t_0,t_1,\bar t_1,t_2,\bar t_2,...)$ be such that
$t_i, \bar t_i=0 $ for $i>n, 0<t_0<1$ and $|t_i|,\ | \bar t_i |
\leqslant (4n^3 2^n e^n)^{-1}$. Then the series $F(\tilde t)$ is convergence.
\endproclaim

A proof is based on a sequential of estimations of all values,
utilized in the definition of $N^2$. Present these estimations:

\it {1. Let $i + j \ = \ s_1 + \dots + s_m$. Then
$P _ {ij} (s_1..., s_m) \leqslant\min (C^{m-1}_{i-1}, C^{m-1}_{j-1})$.}

\demo {Proof}
$$P_{j, i}(s_1, \dots, s_m) \ =P_{i, j}(s_1, \dots, s_m) \ 
= \ \# \{ (i_1, \dots, i_m) \ | \ i \ =$$
$$=\ i_1 + \dots + i_m,\ 1 \leqslant i_r \leqslant s_r-1 \} \ \leqslant$$
$$\leqslant \ \# \{ (i_1, \dots, i_m) | i =
i_1 + \dots + i_m, 1 \leqslant i_r \} \ = \ C_{i-1}^{m-1}\ . \ \qed $$
\enddemo

\it {2. Let $i + j \ = \ s_1 + \dots + s_m$. Then
$T^1_{ij}(s_1,...,s_m)\leqslant\frac{\ell^{m-1}}{m!},\
\text {where} \ \ell =\min (i, j) $.}

\demo{Proof}
$$T_{j, i}^1 (s_1, \dots, s_m) \ =T_{i, j}^1 (s_1, \dots, s_m) \
=\ \sum\limits_{\Sb k \geqslant 1 \\ n_1 + \dots + n_k = m \\ n_r
\geqslant 1 \endSb}\frac {1}{k  n_1!  \dots  n_k!}\ \times$$
$$\times \ P_{i, j}\left(\underbrace{s_1 + \dots + s_{n_1}}_{n_1},
\dots, \underbrace{s_{n_1 + \dots + n_{k-1} + 1} + \dots + s_{n_1 +
\dots + n_k}}_{n_k}\right) \ \leqslant$$
$$\leqslant \ \sum\limits_{\Sb k \geqslant 1 \\ n_1
+ \dots + n_k = m \\ n_r \geqslant 1
\endSb} \frac {C_{i-1}^{k-1}}{k  n_1!  \dots  n_k!} \ =$$
$$= \ \frac{1}{m!}\ \frac{\partial^m}{\partial x^m}\
\frac{1}{i}  \sum\limits_{k \geqslant 1} C_{i}^k  \sum
\limits_{n_1 + \dots + n_k = m} \frac {1}{n_1! \dots n_k!} \ =\ $$
$$\ \frac{1}{m!}\ \frac{\partial^m}{\partial x^m}\ \frac{1}{i}
\sum\limits_{k \geqslant 1}  C_{i}^k  \left.
\left(x + \frac{x^2}{2!} + \dots \right)^k \right|_{x=0} \ =$$
$$=\ \frac{1}{m!}\ \frac{\partial^m}{\partial x^m}\
\frac{1}{i}  \left. \left( \sum\limits_{k \geqslant 1}
C_{i}^k (e^x-1)^k \right) \right|_{x=0} \ =$$
$$= \ \frac{1}{m!}\ \frac{\partial^m}{\partial x^m}\
\frac{1}{i}  \left. ((1+(e^x-1))^i-1) \right|_{x=0} \
= \ \frac{1}{m!}\frac{\partial^m}{\partial x^m}\
\frac{1}{i}  \left. (e^{ix}-1) \right|_{x=0} \ =
\ \frac{i^{m-1}}{m!}\ . \ \qed$$
\enddemo

\it {3.  Let $\ i_1 + \dots + i_k\ = \ s_1 + \dots + s_m$
and $ (\ell_1-1) + \dotsb + (\ell_m-1) =k-2$ . Then
$T^2 _ {i_1..., i_k} \pmatrix s_1 \dotsb s_m \\ \ell_1
\dotsb\ell_m\endpmatrix\ \leqslant\frac{I^{m-1}(k-1)^m(k-2)!}{m!}\ ,$ 
where $I =\max (i_r) $.}

\demo {Proof}
We use an induction by $k$. If $k=2$, then
$$T_{i_1, i_2}^2 \left(\matrix s_1,&\dots,&s_m \\ 1,&\dots,&1
\endmatrix\right) \ = \
T_{i_1, i_2}^1 (s_1, \dots, s_m) \ \leqslant$$
$$\leqslant \ \frac{I^{m-1}}{m!} \
= \ \left.\frac{I^{m-1}(k-1)^{m}(k-2)!}{m!} \right|_{k \ =\ 2}\ .$$

Let $k>2$. Note at first that, if $l = (l_i-1) + \dots + (l_j-1)$, then 
\nopagebreak
$$\sum\limits_{1 \leqslant i \leqslant j \leqslant m,\ j-i=d}
\frac{l}{(k-2)(d+1)} \ = \ \frac{\sum\limits_{1 \leqslant i
\leqslant j \leqslant m,\ j-i=d}  (l_i-1)  +
\dots  +  (l_j-1)}{(k-2)(d+1)} \ \leqslant$$
$$\leqslant \ \frac{\sum\limits_{1 \leqslant t \leqslant m}
(l_t-1) (d+1)}{(k-2)(d+1)} \ = \ 1.$$
Thus
\nopagebreak
$$T_{i_1, \dots, i_k}^2 \left(\matrix s_1,&\dots,&s_m \\
l_1,&\dots,&l_m \endmatrix\right) \ =
\ \sum\limits_{\Sb 1 \leqslant i \leqslant j \leqslant m \\ s,\ l
\geqslant 1 \endSb} l \ T_{s, i_k}^1 (s_i, \dots, s_j) \ \times$$
$$\times \ T_{i_1, \dots, i_{k-1}}^2 \left(\matrix
s_1,&\dots,&s_{i-1},&s,&s_{j+1},&\dots,&s_m \\ l_1,&
\dots,&l_{i-1},&l,&l_{j+1},&\dots,&l_m \endmatrix\right) \ \leqslant$$

$$\leqslant \ \sum\limits_{1 \leqslant i \leqslant j \leqslant m} l
\frac{I^{j-i}}{(j-i+1)!}  \frac{I^{m-(j-i+1)+1-1} 
(k-2)^{m-(j-i+1)+1}  (k-3)!}{(m-(j-i+1)+1)!} \ =$$
$$= \ I^{m-1}  (k-3)!  \sum\limits_{1 \leqslant i \leqslant j \leqslant m} l
 \frac{(k-2)^{m-(j-i)}}{(j-i+1)!  (m-(j-i))!} \ = $$
$$= \ \frac{I^{m-1}  (k-2)!}{m!}  \sum\limits_{1 \leqslant i
\leqslant j \leqslant m,\ d=j-i}  \frac{l}{(k-2)(d+1)}  \frac{m!}{d!
 (m-d)!}  (k-2)^{m-d} \ =$$
$$= \ \frac{I^{m-1}  (k-2)!}{m!}  \sum\limits_{d=0}^{m-1}
\frac{m!}{d!  (m-d)!}  (k-2)^{m-d}  \left(
\sum\limits_{1 \leqslant i \leqslant j \leqslant m,\ j-i=d} 
\frac{l}{(k-2)(d+1)} \right) \ \leqslant$$
\nopagebreak
$$\leqslant \ \frac{I^{m-1}  (k-2)!}{m!}  \sum\limits_{d=0}^{m-1}
\frac{m!}{d!  (m-d)!}  (k-2)^{m-d} \ \leqslant
\ \frac{I^{m-1}  (k-2)!}{m!}  \sum\limits_{t=0, t=m-d}^{m}
C_m^{t}  (k-2)^{t} \ =$$
$$= \ \frac{I^{m-1}  (k-2)!}{m!}  ((k-2)+1)^m \ =
\ \frac{I^{m-1}  (k-1)^m  (k-2)!}{m!}\ . \ \qed$$
\enddemo

\it {4. Let $\bar I = max_r (\bar i_r)$ and 
$$\tilde S_{\bar i_1, \dots, \bar i_{\bar k}} (m, k) \ =
\sum\limits_{\Sb \{\bar i_1^1, \dots, \bar i_1^{n_1}\}
\sqcup \dots \sqcup \{\bar i_m^1, \dots, \bar i_m^{n_m}\} =\\
= \{\bar i_1, \dots, \bar i_{\bar k}\} \\ (l_1-1) + \dots + (l_m-1)
= k-2 \\ n_r,\ l_r \geqslant 1 \\ s_r=\bar i_r^1 + \dots + \bar i_r^{n_r}
\\ s_r-n_r-\ell_r+1\geqslant 0 \endSb}
\frac{(s_1-1)!}{(s_1-n_1-l_1+1)!  (l_1-1)!} \times \dots \times$$
$$\times \ \frac{(s_m-1)!}{(s_m-n_m-l_m+1)!  (l_m-1)!}\ .$$

Then $\tilde S_{\bar i_1, \dots, \bar i_{\bar k}}
(m, k) \ \leqslant \ m  (\bar k-1)!  C_{\bar I \bar k - \bar k}^{k-2}
 C_{\bar I \bar k}^{\bar k-m}$.}

\demo{Proof} We use the equality
$$\sum\limits_{\tilde n_1 + \dots + \tilde n_m = \bar k-m}
C_{\bar I \tilde n_1 + \bar I}^{\tilde n_1} \times \dots \times
C_{\bar I \tilde n_m + \bar I}^{\tilde n_m}  \frac{\bar I}
{\bar I \tilde n_1 + \bar I} \times \dots \times \frac{\bar I}
{\bar I \tilde n_m + \bar I} \ =$$
$$=\ C_{\bar I (\bar k-m) + m \bar I}^{\bar k-m} 
\frac{m\bar I}{\bar I (\bar k-m) + m\bar I}$$
that follows from [8,5.63].
Then
$$\tilde S_{\bar i_1, \dots, \bar i_{\bar k}} (m, k) \ =$$
$$= \ \sum\limits_{\Sb \{\bar i_1^1, \dots, \bar i_1^{n_1}\}
\sqcup \dots \sqcup \{\bar i_m^1, \dots, \bar i_m^{n_m}\}
=\\= \{\bar i_1, \dots, \bar i_{\bar k}\} \\
(l_1-1) + \dots + (l_m-1) = k-2 \\ n_r,\ l_r \geqslant 1 \\ s_r
= \bar i_r^1 + \dots + \bar i_r^{n_r}
\\ s_r-n_r-\ell_r+1\geqslant 0 \endSb}
\frac{(s_1-1)!}{(s_1-n_1-l_1+1)!  (l_1-1)!} \times \dots \times $$
$$\times \frac{(s_m-1)!}{(s_m-n_m-l_m+1)!  (l_m-1)!} \ =$$
$$= \ \sum\limits_{\Sb \{\bar i_1^1, \dots, \bar i_1^{n_1}\}
\sqcup \dots \sqcup \{\bar i_m^1, \dots, \bar i_m^{n_m}\}
=\\= \{\bar i_1, \dots, \bar i_{\bar k}\} \\
s_r = \bar i_r^1 + \dots + \bar i_r^{n_r}
\geqslant n_r \geqslant 1 \endSb}
\frac{(s_1-1)!}{(s_1-n_1)!} \times \dots \times
\frac{(s_m-1)!}{(s_m-n_m)!}  \times$$
$$\times  \sum\limits_{\Sb (l_1-1) + \dots + (l_m-1) = k-2
\\ s_r-n_r-\ell_r+1\geqslant 0 \endSb}
\frac{(s_1-n_1)!}{(s_1-n_1-l_1+1)!  (l_1-1)!} \times \dots\times$$
$$\times\ \frac{(s_m-n_m)!}{(s_m-n_m-l_m+1)!  (l_m-1)!} \ =$$
$$= \ \sum\limits_{\Sb \{\bar i_1^1, \dots, \bar i_1^{n_1}\}
\sqcup \dots \sqcup \{\bar i_m^1, \dots, \bar i_m^{n_m}\}
=\\= \{\bar i_1, \dots, \bar i_{\bar k}\}
\\ s_r = \bar i_r^1 + \dots + \bar i_r^{n_r}
\geqslant n_r \geqslant 1 \endSb}
\frac{(s_1-1)!}{(s_1-n_1)!} \times \dots \times$$
$$\times \frac{(s_m-1)!}{(s_m-n_m)!}
\sum\limits_{\tilde l_1 + \dots + \tilde l_m = k-2}
C_{s_1-n_1}^{\tilde l_1}  \dots  C_{s_m-n_m}^{\tilde l_m} \ =$$
$$= \ \sum\limits_{\Sb \{\bar i_1^1, \dots, \bar i_1^{n_1}\}
\sqcup \dots \sqcup \{\bar i_m^1, \dots, \bar i_m^{n_m}\}
=\\= \{\bar i_1, \dots, \bar i_{\bar k}\}
\\ s_r = \bar i_r^1 + \dots + \bar i_r^{n_r}
\geqslant n_r \geqslant 1  \endSb}
\frac{(s_1-1)!}{(s_1-n_1)!} \times \dots \times$$
$$\times\ \frac{(s_m-1)!}{(s_m-n_m)!}  C_{(s_1+\dots+s_m)-
(n_1+\dots+n_m)}^{k-2} \ =$$
$$= \ C_{\bar i_1 + \dots + \bar i_{\bar k} - \bar k}^{k-2} 
\sum\limits_{\Sb \{\bar i_1^1, \dots, \bar i_1^{n_1}\}
\sqcup \dots \sqcup \{\bar i_m^1, \dots, \bar i_m^{n_m}\}
=\\= \{\bar i_1, \dots, \bar i_{\bar k}\}
\\ s_r = \bar i_r^1 + \dots + \bar i_r^{n_r}
\geqslant n_r \geqslant 1 \endSb}
\frac{(s_1-1)!}{(s_1-n_1)!} \times \dots \times
\frac{(s_m-1)!}{(s_m-n_m)!}\ \leqslant\ $$
$$\leqslant\ \tilde S_{\underbrace{\bar I, \dots, \bar I}_{\bar k}} (m, k) \
= \ C_{\bar I \bar k - \bar k}^{k-2}  \sum\limits_{\Sb
\{\bar i_1^1, \dots, \bar i_1^{n_1}\} \sqcup \dots \sqcup
\\ \sqcup \{\bar i_m^1, \dots, \bar i_m^{n_m}\}
=\\= \{\bar I, \dots, \bar I\},\ n_r \geqslant 1 \\ s_r
= \bar i_r^1 + \dots + \bar i_r^{n_r} = I n_r\endSb}
\frac{(s_1-1)!}{(s_1-n_1)!} \times \dots \times
\frac{(s_m-1)!}{(s_m-n_m)!} \ =$$
$$= \ C_{\bar I \bar k - \bar k}^{k-2}  \sum\limits_{\Sb
n_1 + \dots + n_m = \bar k \\ n_r \geqslant 1 \endSb}
\frac{\bar k!}{n_1!  \dots  n_m!} 
\frac{(\bar I n_1 - 1)!}{(\bar I n_1 - n_1)!} \times \dots \times
\frac{(\bar I n_m - 1)!}{(\bar I n_m - n_m)!} \ =$$
$$= \ \bar k!  C_{\bar I \bar k - \bar k}^{k-2} 
\sum\limits_{\Sb \tilde n_1 + \dots + \tilde n_m = \bar k - m
\\ \tilde n_r = n_r -1 \geqslant 0 \endSb}
\frac{(\bar I \tilde n_1 + \bar I - 1)!}{(\tilde n_1+1)!
(\bar I \tilde n_1 + \bar I - \tilde n_1 -1)!} \times \dots\times$$
$$\times\ \frac{(\bar I \tilde n_m + \bar I - 1)!}{(\tilde n_m+1)!
 (\bar I \tilde n_m + \bar I - \tilde n_m -1)!} \ =$$
$$= \ \bar k!  C_{\bar I \bar k - \bar k}^{k-2} 
\sum\limits_{\tilde n_1 + \dots + \tilde n_m =
\bar k - m} \frac{\bar I \tilde n_1 + \bar I -
\tilde n_1}{(\bar I \tilde n_1 + \bar I)  (\tilde n_1+1)}
 \frac{(\bar I \tilde n_1 + \bar I)!}{\tilde n_1! 
(\bar I \tilde n_1 + \bar I - \tilde n_1)!} \times \dots \times$$
$$\times\ \frac{\bar I \tilde n_m + \bar I - \tilde n_m}{(\bar I
\tilde n_m + \bar I)  (\tilde n_m+1)}  \frac{(\bar I
\tilde n_m + \bar I)!}{\tilde n_m!  (\bar I \tilde n_m +
\bar I - \tilde n_m)!} \ \leqslant$$
$$\leqslant \ \bar k!  C_{\bar I \bar k - \bar k}^{k-2} 
\sum\limits_{\tilde n_1 + \dots + \tilde n_m = \bar k - m}
\frac{\bar I}{\bar I \tilde n_1 + \bar I}  C_{\bar I \tilde n_1
+ \bar I}^{\tilde n_1} \times \dots \times
\frac{\bar I}{\bar I \tilde n_m + \bar I}  C_{\bar I \tilde n_m
+ \bar I}^{\tilde n_m} \ =$$
$$= \ \bar k!  C_{\bar I \bar k - \bar k}^{k-2} 
C_{\bar I (\bar k-m) + m \bar I}^{\bar k-m} 
\frac{m\bar I}{\bar I (\bar k-m) + m\bar I} \ \leqslant
\ m  (\bar k-1)!  C_{\bar I \bar k - \bar k}^{k-2} 
C_{\bar I \bar k}^{\bar k-m}\ . \ \qed$$
\enddemo

\it {5. $N^1_i (i_1..., i_k |\bar i_1..., \bar i _ {\bar k}) \ \leqslant
\ (k-1)! (\bar k-1)! e^{I(k-1)} 2^{\bar I \bar k - \bar k}  2^{\bar I
\bar k}.$}

\demo {Proof}
$$N_i^1 (i_1, \dots, i_k  |  \bar i_1, \dots, \bar i_{\bar k}) \ = $$
$$= \ \sum\limits_{\Sb m \geqslant 1 \\ s_1 + \dots + s_m = i_1 + \dots +
i_k \\ l_1 + \dots + l_m = m+k-2 \\ s_r, l_r \geqslant 1 \endSb}
(-1)^{m+1} \ S_{\bar i_1, \dots, \bar i_{\bar k}} \left(
\matrix s_1,&\dots,&s_m \\ l_1,&\dots,&l_m \endmatrix\right)  \times$$
$$\times \ T_{i_1, \dots, i_k}^2 \left(\matrix
s_1,&\dots,&s_m \\ l_1,&\dots,&l_m \endmatrix\right) \ \leqslant$$
$$\leqslant\ \sum\limits_{\Sb m \geqslant 1\\s_1
+\dots+s_m=i_1+\dots+i_k\\l_1+\dots+l_m
=m+k-2 \\ s_r,l_r \geqslant 1 \endSb}
\frac{I^{m-1} (k-1)^{m}  (k-2)!}{m!}  \times$$
$$\times  \sum\limits_{\Sb \{\bar i_1^1, \dots,
\bar i_1^{n_1}\} \sqcup \dots \sqcup \{\bar i_m^1, \dots,
\bar i_m^{n_m}\} =\\= \{\bar i_1, \dots, \bar i_{\bar k}\} \\
\bar i_r^1 + \dots + \bar i_r^{n_r} = s_r
\\s_r-n_r-\ell_r+1\geqslant 0 \endSb}
\frac{(s_1-1)!}{(s_1-n_1-l_1+1)!  (l_1-1)!} \times \dots \times$$
$$\times \ \frac{(s_m-1)!}{(s_m-n_m-l_m+1)!  (l_m-1)!} \ =$$
$$= \ \sum\limits_{m \geqslant 1}
\frac{I^{m-1} (k-1)^{m-1} (k-2)!}{m!}  \times$$
$$\times  \sum\limits_{\Sb \{\bar i_1^1, \dots,
\bar i_1^{n_1}\} \sqcup \dots \sqcup \{\bar i_m^1,
\dots, \bar i_m^{n_m}\} =\\= \{\bar i_1, \dots,
\bar i_{\bar k}\} \\ (l_1-1) + \dots + (l_m-1) =
k-2 \\ n_r,\ l_r \geqslant 1 \\ s_r = \bar i_r^1 + \dots
+ \bar i_r^{n_r} \\s_r-n_r-\ell_r+1\geqslant 0 \endSb}
\frac{(s_1-1)!}{(s_1-n_1-l_1+1)! (l_1-1)!} \times \dots\times$$
$$\times\frac{(s_m-1)!}{(s_m-n_m-l_m+1)! (l_m-1)!} \ =$$
$$= \ \sum\limits_{m \geqslant 1}
\frac{I^{m-1} (k-1)^{m-1} (k-2)!}{m!} \ \tilde S_{\bar i_1,
\dots, \bar i_{\bar k}} (m, k)\ \leqslant$$
$$\leqslant \ \sum\limits_{m \geqslant 1}
\frac{I^{m-1} (k-1)^m (k-2)!}{m!} \
\tilde S_{\bar i_1, \dots, \bar i_{\bar k}} (m, k) \ \leqslant$$
$$\leqslant \ \sum\limits_{m \geqslant 1}\frac{I^{m-1} (k-1)^m (k-2)!}{m!}
m (\bar k-1)!  C_{\bar I \bar k - \bar k}^{k-2}
C_{\bar I \bar k}^{\bar k-m} \ =$$
$$=\ (k-1)! (\bar k-1)! \ \sum\limits_{m \geqslant 1}\frac{I^{m-1}
(k-1)^{m-1}}{(m-1)!}\ C_{\bar I \bar k - \bar k}^{k-2}
\ C_{\bar I \bar k}^{\bar k-m}\ \leqslant$$
$$\leqslant \ (k-1)!  (\bar k-1)!  e^{I(k-1)}  2^{\bar I \bar k - \bar k}
2^{\bar I \bar k}\ . \qed$$
\enddemo

\demo{Proof of theorem 3}
The coefficient for
$t_0^{\dots} t_{i_1}^{n_1} \dots t_{i_I}^{n_I} \bar
t_{\bar i_1}^{\bar n_1} \dots \bar t_{\bar i_{\bar I}}^{\bar n_{\bar I}}$
is equal
$$\frac{i_1^{n_1}  \dots  i_I^{n_I}}{n_1!  \dots  n_I!} \
\frac{\bar i_1^{\bar n_1}  \dots  \bar i_{\bar I}^{\bar n_{\bar I}}}
{\bar n_1!  \dots  \bar n_{\bar I}!}
\ N_i^2 \left(\matrix i_1,&\dots,&i_I \\ n_1,&\dots,&n_I \endmatrix \left|
\matrix \bar i_1,&\dots,&\bar i_{\bar I} \\ \bar n_1,&\dots,
&\bar n_{\bar I} \endmatrix \right.\right) \ =$$
$$= \ \frac{i_1^{n_1}  \dots  i_I^{n_I}}{n_1!  \dots 
n_I!} \ \frac{\bar i_1^{\bar n_1}  \dots 
\bar i_{\bar I}^{\bar n_{\bar I}}}{\bar n_1!  \dots  \bar n_{\bar I}!}
\ N_i^1 \left(\underbrace{i_1, \dots, i_1}_{n_1}, \dots,
\underbrace{i_I, \dots, i_I}_{n_I} \left|
\underbrace{\bar i_1, \dots, \bar i_1}_{\bar n_1}, \dots,
\underbrace{\bar i_{\bar k}, \dots,
\bar i_{\bar I}}_{\bar n_{\bar I}} \right.\right)\ \leqslant$$

$$\leqslant \ \frac{i_1^{n_1}  \dots  i_I^{n_I}}{n_1!  \dots  n_I!} \
\frac{\bar i_1^{\bar n_1}  \dots 
\bar i_{\bar I}^{\bar n_{\bar I}}}{\bar n_1!  \dots 
\bar n_{\bar I}!} 
k!  \bar k!  e^{\tilde I (k-1)}  2^{2\tilde I \bar k} \ \leqslant
\ \tilde I^K  e^{\tilde I K}  2^{\tilde I K}  \frac{k!  \bar k!}{n_1! 
\dots  n_I!  \bar n_1!  \dots  \bar n_{\bar I}!}\ \leqslant $$
$$\ \leqslant \ \tilde I^K  e^{\tilde I K}  2^{\tilde I K}
\tilde I^K \ \leqslant \ (\tilde I^2 2^{\tilde I} e^{\tilde I})^K\ , $$
where
$k=n_1+\dots +n_I$, $\bar k=\bar n_1+\dots +\bar n_I$, $K=k+\bar k$ and 
$\tilde I= max(I,\bar I)$.

Consider now monomials from $t_0,t_1,\bar t_1,\dots ,t_n,\bar t_n$ of degree
$K$. The number of such monomial is $(2n)^K$. Thus its sum in the series is
not more that 
$$(n^2  2^n  e^n)^K  (2n)^K  (4  n^3  2^n e^n)^{-K} \ \leqslant 2^{-K}\ .$$
This implies the convergence of the series $F(\tilde t)$.
\qed
\enddemo

\end